\documentclass[aps,prd,preprint2,groupedaddress]{revtex4-1}

\usepackage{natbib}

\usepackage{amsmath}
\usepackage{amssymb}
\usepackage{graphicx}

\newcommand{\bkap}{{\boldsymbol a}}

\newcommand{\nhat}{{\hat {\bf n}}}
\newcommand{\be}{{\hat {\bf e}}}
\newcommand{\kapxx}{a_{xx}}
\newcommand{\kapyy}{a_{yy}}
\newcommand{\kapxy}{a_{xy}}
\newcommand{\kapyx}{a_{yx}}
\newcommand{\bnabla}{{\boldsymbol \nabla}}
\newcommand{\bx}{{\boldsymbol x}}
\newcommand{\bk}{{\boldsymbol k}}

\begin{document}

\title{Renormalization Group theory outperforms other approaches in statistical comparison between upscaling techniques for porous media} 

\author{Shravan Hanasoge}
\email[]{hanasoge@tifr.res.in}
\affiliation{Department of Astronomy and Astrophysics, Tata Institute of Fundamental Research, Mumbai 400005, India}
\author{Umang Agarwal and Kunj Tandon}
\affiliation{Shell India Markets Pvt.~Ltd.,~Bangalore Hardware Park, Devanahalli Industrial Park,  Bengaluru 562149, India}
\author{J. M. Vianney A. Koelman}
\affiliation{Center for Computational Energy Research, Eindhoven University of Technology / Dutch Institute For Fundamental Energy Research, Eindhoven, The Netherlands}

\date{\today}

\begin{abstract}
Determining the pressure differential required to achieve a desired flow rate in a porous medium requires solving Darcy's law, a Laplace-like equation, with a spatially varying tensor permeability. In various scenarios, the permeability coefficient is sampled at high spatial resolution, which makes solving Darcy's equation numerically prohibitively expensive. As a consequence, much effort has gone into creating upscaled or low-resolution effective models of the coefficient while ensuring that the estimated flow rate is well reproduced, bringing to fore the classic tradeoff between computational cost and numerical accuracy. Here we perform a statistical study to characterize the relative success of upscaling methods on a large sample of permeability coefficients that are above the percolation threshold. We introduce a new technique based on Mode-Elimination Renormalization-Group theory (MG) to build coarse-scale permeability coefficients. Comparing the results with coefficients upscaled using other methods, we find that MG is consistently more accurate, particularly so due to its ability to address the tensorial nature of the coefficients. MG places a low computational demand, in the manner that we have implemented it, and accurate flow-rate estimates are obtained when using MG-upscaled permeabilities that approach or are beyond the percolation threshold.
\end{abstract}


\maketitle

\section{Introduction}\label{intro}
Hydrocarbons are often recovered from reservoirs by applying a pressure differential across a pair of locations and extracting the output fluid from one. The reservoir performance is dependent on the flow behaviour of hydrocarbon and pumped fluid through the reservoir section. The heterogeneities across different length scales in reservoir formations is typically captured in the reservoir simulations by an effective gridblock value. This effective value is derived using an upscaling process where the heterogeneous medium represented by the upscaled grid block would essentially have the same flow at its boundaries under the same pressure graident as the finer scale medium. The physics of fluid flow through a porous medium, described by the time-independent Darcy's law is,
\begin{equation}
\bnabla\cdot(\bkap\cdot\bnabla\phi) = 0,\label{darcy}
\end{equation}
where $\phi = \phi(\bx)$ is the pressure distribution, $\bx$ is the spatial coordinate, $\bkap = \kapxx \be_x \be_x +  \kapxy  \be_x \be_y +  \kapyx  \be_y \be_x + \kapyy  \be_y \be_y$ is the tensor permeability, $\be_x$ and $\be_y$ are unit vectors along the $x$ and $y$ directions and $\bnabla$ is the spatial gradient. Although we limit the current discussion to 2-D, equation~(\ref{darcy}) applies in 3-D as well. The specific quantity of interest to us is the flow rate, defined by
\begin{equation}
f = -\int_S dA\, \nhat\cdot \bkap\cdot\bnabla\phi,
\end{equation}
where $S$ represents the surface perpendicular to the direction in which the pressure differential is applied, $\nhat$ is the normal vector to that surface (and therefore co-aligned with the direction of the pressure differential) and the area integral is over the entire surface $S$. The negative sign indicates that the flow is in the direction of decreasing pressure. The solution to equation~(\ref{darcy}) therefore gives the flow rate $f$ associated with a specific pressure differential. Since the equation is linear, the solution $\phi$ scales directly with the magnitude of the pressure difference - i.e., the flow rate obtained through the solution of equation~(\ref{darcy}) is directly proportional to the pressure contrast.  

Contemporary methodologies  \citep{king89,renard1997,yeo2001,lunati2001,farmer2002,green2007,Khalili2012,Li2016} allow for building models of the spatial distribution of permeabilities in reservoirs at high spatial resolution, resulting in $\sim1024^3$-sized blocks in 3D. Additionally, the permeabilities may display large contrasts and both issues render the computation intractable (or difficult, at any rate), and certainly so on a routine basis. The commonly adopted approach is to design an effective coarse version of the coefficient such that low-wavenumber aspects of the solution are preserved, thereby ensuring that the flow rate changes minimally.   A variety of prevailing approaches are applied to derive the coarsened or upscaled coefficient from the original high-resolution representation. The simplest method is to directly filter out high wavenumbers from the coefficient and project it on to a coarser grid. Depending on the types of spatial variations present in the coefficient, this may be an accurate approximation but success is by no means always assured. Indeed, for the types of spatial variations seen in realistic coefficients, straightforward averaging may not be the best technique of upscaling \citep[e.g.][]{king89}. Additionally, simple averaging does not explicitly account for how the spatial geometry and tensor components of the coefficient participate in influencing the solution. {In particular, the solution $\phi$ to equation~(\ref{darcy}) is sensitively dependent on both diagonal ($\kapxx$ and $\kapyy$) and off-diagonal components ($\kapxy$). Most upscaling methodologies do not set out a  procedure to accommodate the off-diagonal component, choosing rather to ignore them, i.e. $\kapxy = 0$. Depending on the magnitude of the off-diagonal relative to the diagonal component, this assumption may lead to potentially significant errors and systematical biases in estimates of flow rates.}

More sophisticated techniques such as Renormalization Group Theory and the technique of \cite{karim}, termed KK hence, explicitly address the tensorial nature of the coefficient (although KK only handles diagonal components of the coefficient tensor). Renormalization Group, which has a celebrated history in particle physics, made possible the construction of effective (low-wavenumber) theories in multi-scale, multi-component systems. Here we use a variant termed Mode-Elimination Renormalization Group Theory (MG), introduced by Kenneth Wilson \citep[and therefore known as the Wilsonian Renormalization Group;][]{wilson1975}. Let us consider the Fourier transform of equation~(\ref{darcy}); the spatial coordinate is replaced by the wavenumber, i.e. $\bx \rightarrow \bk$ and functions of real space by their transforms. In order to minimise notational burden, we use the same symbol to denote the transformed and original quantities, e.g. $\phi(\bk)$ and $\phi(\bx)$ are transform pairs, as are $\bkap(\bk)$ and $\bkap(\bx)$. Consider the discrete Fourier transform of equation~(\ref{darcy}),
\begin{equation}
\sum_{\bx} e^{-i\bk\cdot\bx} \bnabla\cdot(\bkap\cdot\bnabla\phi) =\sum_{\bk'} i\bk\cdot\bkap(\bk-\bk')\cdot i\bk'\,\phi(\bk') = -[{\rm BC}] = -\sum_{\bk'} \bk\cdot\bkap(\bk-\bk')\cdot \bk'\,\phi(\bk') ,\label{darcy.ft0}
\end{equation}
where the spatial gradient Fourier transforms in to the wavenumber, i.e. ${\bnabla} \rightarrow -i\bk$, the spatial product turns into a convolution, and the boundary-condition term (denoted BC) is moved to the right side (thereby acting as a source). Equation~(\ref{darcy.ft0}) shows that the wavenumbers of $\bkap$ and $\phi$ are fully mixed as a consequence of the convolution. Given $N$ discrete wavenumbers, equation~(\ref{darcy.ft0}) may be written as an $N\times N$ matrix inversion problem, $\{ \bk\cdot\bkap(\bk-\bk')\cdot \bk' \} [\phi(\bk')] = [\rm BC]$, where the first term in brackets indicates an $N\times N$ matrix with rows being associated with wavenumbers $\bk$ and columns with $\bk'$ and terms in square brackets denote $N\times1$ column vectors. 
Since we are only interested in ensuring that the low wavenumber components of $\phi$ are faithfully reproduced,
the challenge is to build an appropriate reduced version of the matrix inversion.   
In particular, the diffusion problem has a long history in applied mathematics as well, where it is known as ``homogenisation", gaining prominence through the work of \cite{kozlov79} and \cite{varadhan82}. \cite{brewster} developed a treatment for finite-wavenumber cutoff in this context.

\section{Renormalization Group\label{RG.theory}}
Here we apply MG, among other methods, to develop a coarse-scale representation of the coefficient $\bkap$. Building on the work of \cite{dorobantu}, \cite{hanasoge.RG} analysed the performance of the method in the context of the diffusion equation. We introduce a cutoff wavenumber, denoted by $\bk_c$ and describe wavenumbers that are greater or less than this by $\bk_>$ (all wavenumbers $\bk > \bk_c$) and $\bk_<$ (all wavenumbers $\bk < \bk_c$). Only the coarse-scale part of the solution $\phi(\bk_<)$ interests us, and we do not wish to compute $\phi(\bk_>)$. Since equation~(\ref{darcy.ft0}) holds for each wavenumber $\bk$, we split it into two equations, one each for low and high wavenumbers,
\begin{eqnarray}
-\sum_{\bk'} \bk_<\cdot\bkap(\bk_<-\bk')\cdot \bk'\,\phi(\bk')  = [{\rm BC}]_<\label{high.darcy.ft1}\\
-\sum_{\bk'} \bk_>\cdot\bkap(\bk_>-\bk')\cdot \bk'\,\phi(\bk') = -[{\rm BC}]_>.\label{darcy.ft1}
\end{eqnarray}
If the boundary condition $[{\rm BC}]$ has no high-wavenumber power, the right hand side of equation~(\ref{high.darcy.ft1}) is zero,
\begin{equation}
\sum_{\bk'} \bk_>\cdot\bkap(\bk_>-\bk')\cdot \bk'\,\phi(\bk') = 0.\label{eq.interm}
\end{equation}
The sum over $\bk'$ in equation~(\ref{eq.interm}) applies to the entire range of wavenumbers and may be split, as in equations~(\ref{high.darcy.ft1}) and~(\ref{darcy.ft1}) into two partial sums, i.e. for low and high wavenumbers
\begin{equation}
\sum_{\bk'_>} \bk_>\cdot\bkap(\bk_>-\bk'_>)\cdot \bk'_>\,\phi(\bk'_>) + \sum_{\bk'_<} \bk_>\cdot \bkap(\bk_>-\bk'_<)\cdot \bk'_<\,\phi(\bk'_<) = 0,
\end{equation}
which is ultimately written in matrix form to assist the analysis,
\begin{equation}
\{\bkap_{>>}\} [\phi]_> = - \{\bkap_{><}\} [\phi]_<,
\end{equation}
where $\{\bkap_{>>}\} = \bk_>\cdot\bkap(\bk_> - \bk'_>)\cdot\bk'_>$ is a square matrix of size $N_>\times N_>$, where $N_>$ is the count of high wavenumbers, and $\{\bkap_{><}\} = \bk_>\cdot\bkap(\bk_> - \bk'_<)\cdot\bk'_<$ is of size $N_<\times N_>$, where $N_<$ is the count of the (desired) low wavenumbers, likely a non-square matrix. The elements of these matrices correspond to the allowed values of $\bk_>, \bk'_>$, and $[\phi]_< = \phi(\bk'_<)$ and $[\phi]_> = \phi(\bk'_>)$ are column vectors. This allows us to express the high-spatial-frequency component of the solution, $\phi(\bk'_>)$ in terms of the low-wavenumber part (that we are interested in computing)
\begin{equation}
 [\phi]_> = - \{\bkap_{>>}\}^{-1}\,\{\bkap_{><}\} [\phi]_<.\label{highwave.comp}
\end{equation}
Now we study the low wavenumber part of equation~(\ref{darcy.ft0}), 
\begin{equation}
\sum_{\bk'} \bk_<\cdot\bkap(\bk_>-\bk')\cdot \bk'\,\phi(\bk') = [{\rm BC}],\label{low.wave}
\end{equation}
where for the sake of notational convenience, $[{\rm BC}]_<$ is written as $[{\rm BC}]$. Expanding equation~(\ref{low.wave}) into low and high wavenumbers again,
\begin{equation}
\sum_{\bk'_>} \bk_<\cdot\bkap(\bk_<-\bk'_>)\cdot \bk'_>\,\phi(\bk'_>) + \sum_{\bk'_<} \bk_<\cdot\bkap(\bk_<-\bk'_<)\cdot \bk'_<\,\phi(\bk'_<) = [{\rm BC}], 
\end{equation}
or in matrix notation,
\begin{equation}
\{\bkap_{<<}\} [\phi]_< + \{\bkap_{<>}\} [\phi]_> =  [{\rm BC}].\label{lowwave.comp}
\end{equation}
Note that the same notations as employed in e.g. equation~(\ref{highwave.comp}) are used here as well. We now substitute the expression for $[\phi]_>$ from equation~(\ref{highwave.comp}) in equation~(\ref{lowwave.comp}) to obtain
\begin{equation}
(\{\bkap_{<<}\}  - \{\bkap_{<>}\} \{\bkap_{>>}\}^{-1}\,\{\bkap_{><}\})[\phi]_< =  [{\rm BC}],\label{renorm}
\end{equation}
which is an effective equation for solely the low wavenumbers of $\phi$, or $[\phi]_<$, as we term it here. Equation~(\ref{renorm}) encapsulates the principles of mode-elimination Renormalization group theory. It represents computing the {\it Schur complement} of the full $N\times N$ matrix $\{\bk\cdot\bkap\cdot\bk'\}$ and is equivalently known as numerical homogenisation (or homogenisation) in applied mathematics.

Up to this point, there are no assumptions, except with regards periodic boundary conditions, which may be relaxed by using a non-Fourier basis. In principle, this should give us the exact low-frequency solution, with the added bonuses that computing this portion of $\phi$ is ostensibly computationally tractable and equation~(\ref{renorm}) fully accounts for the tensor nature of $\bkap$. 
The universe hands out no free lunches and MG is no different, firstly for the reason that computing the Schur complement is very expensive, specifically when attempting to substantially reduce the size of the matrix. Secondly, owing to an uncertainty principle that governs convolutional operators \citep{hanasoge.RG}, what was once a pure differential equation turns into a mixed integro-differential equivalent, and therefore difficult to solve. 

As a consequence, we adopt an empirical approach broadly based on MG, where we gradually decimate (reduce in size) the coefficient matrix into a coarse version. The high-resolution coefficient is divided up into blocks of a given size, and each block is decimated into one coefficient and the process is repeated. The size of the block, the eventual size of the coefficient and the decimation technique are choices that must be determined. 

\section{Problem setup}
To develop insight into the performance of upscaling techniques, we limit the problem to 2D. We operate under the assumption that we learn in 2D can be naturally extended to 3D. The first step lies in developing a means to obtain solutions to equation~(\ref{darcy}), which we solve in real space. 
\subsection{Solving Darcy's equation in 2D}
The conditions we apply are $\phi(x=0, y) = 1$, $\phi(x=1, y) = 0$ on the horizontal sides and $\partial\phi/\partial y |_ {(x,y=0)} = 0 = \partial\phi/\partial y |_ {(x,y=1)}$ on the vertical boundaries. We expand equation~(\ref{darcy}) as
\begin{equation}
(\partial_x\kapxx + \partial_y\kapyx)\,\frac{\partial\phi}{\partial x} + (\partial_x\kapxy + \partial_y\kapyy)\,\frac{\partial\phi}{\partial y} +  \kapxx\frac{\partial^2\phi}{\partial x^2} +  \kapyy\frac{\partial^2\phi}{\partial y^2} +  2\kapxy\frac{\partial^2\phi}{\partial x \partial y} = 0,\label{full.eq}
\end{equation}
where the tensor components of the coefficient are $\bkap = \{\kapxx, \kapxy, \kapyx, \kapyy\}$ and $\kapxy = \kapyx$ is honoured. 
We use the following centred fourth-order stencils to resolve first- and second-order derivatives in the $x$ and $y$ directions although only shown here for the $x$ derivative,
\begin{eqnarray}
\frac{\partial\phi}{\partial x}^{i,j} = a_2 (\phi^{i+2 j} - \phi^{i-2 j}) + a_1(\phi^{i+1 j} - \phi^{i-1 j}),\\
\frac{\partial^2\phi}{\partial x^2}^{i,j} = b_2 (\phi^{i+2 j} + \phi^{i-2 j}) + b_1(\phi^{i+1 j} + \phi^{i-1 j}) + b_0\, \phi^{ij},
\end{eqnarray}
where $a_1 = 8/12, a_2 = -1/12$ and $b_0= -30/12 ,b_1 = 16/12, b_2 = -1/12$ and $i,j$ denote $x$ and $y$ grid-point indices respectively. Because the normal derivative $\partial\phi/\partial y = 0$ on the upper and lower boundaries, as the conditions dictate, we introduce ghost nodes. For instance, on the lower boundary, this would imply $\phi^{i, -1} = \phi^{i,1}$ and $\phi^{i, -2} = \phi^{i,2}$. Because we apply Dirichlet conditions on the horizontal sides, we do not introduce ghost nodes; instead, the derivative drops to second order at points adjacent to the $x$ boundaries. The $x$-boundary conditions themselves are shifted to the right-hand side. Standard centred second-order accurate schemes are used for the first and second derivatives (we do not list them here).

With this approach, we reformulate the differential equation as a classic matrix-inversion problem of the sort $A[\phi] = b$, where in this case, $A$ is a banded matrix comprising the discretised terms of equation~(\ref{full.eq}) and $b$ is the boundary condition. Note that the bandwidth of $A$ directly scales with the order of accuracy of the derivative scheme. Once the matrix is constructed, we use the Portable Extensible Scientific Toolkit \citep[PETSc,][]{petsc-main}, a freely available library to solve the matrix equation. PETSc provides a variety of preconditioners such as geometric and algebraic multigrid, and numerous options for iterative Krylov solvers. 
Because we are dealing with a sufficiently small problem in 2D, and given that iterative solutions take a long time to converge, we 
rely on an LU decomposition algorithm to directly solve for $\phi$. Currently PETSc solves $512\times512$ sized problems with high accuracy quickly, but extending to larger blocks possibly requires preconditioned iterative solves (and is an avenue for future work). We work here with the largest possible grid size, i.e. $512\times512$.

\subsection{Solution validation and quantification of error}
Since the geometry of the problem is such that the flow is in the $x$ direction, the flow rate $f$ is given by
\begin{equation}
f = -\int_0^1 dy\,\left( \kapxx \frac{\partial\phi}{\partial x} +  \kapxy \frac{\partial\phi}{\partial y} \right).\label{f.x}
\end{equation}
Because the boundary condition on the left side is $\phi = 1$ and on the right side is $\phi = 0$, $\partial_x \phi < 0$ over most of the domain, which results in a negative value of the integral and therefore a positive value for $f$. 

To ensure that the solution is numerically accurate, we may compute the flow rate $f$ from equation~(\ref{f.x}), seen to be a function of the $x$ coordinate since it is an integral over the entire $y$ range. However the flow rate must be an invariant quantity and therefore $f$ has to be constant. To validate the solution, we calculate $f^{\rm max} / f^{\rm min}$, i.e.
the maximum-to-minimum ratio of the flow rate $f$ as obtained from estimates at all the $x$ coordinates. If this ratio is smaller than some critical threshold (user prescribed), we declare the corresponding solution admissible. The solution associated with the full coefficients (i.e. with no upscaling) is more susceptible to numerical errors because of the presence of sharper and spatially complicated features in the permeability coefficients. For a given coefficient, we therefore only require the exact solution to pass this test.

Given these solutions, the primary question of interest is: do the upscaled coefficient models produce the right flow rate? If not, how erroneous are the estimates? We define the error $\varepsilon$ as
\begin{equation}
\varepsilon = \frac{f^{\rm exact} - f^{\rm model}}{f^{\rm exact}} \times 100,\label{error.def}
\end{equation}
allowing us to characterise the performance of upscaling methods.

\subsection{Generating permeabilities}
In realistic scenarios, the coefficient $\bkap$ may display large contrasts, and some regions in a reservoir are altogether impermeable. Setting $\bkap = {\bf 0}$ or to an extremely small value in a numerical calculation signals trouble because the matrix becomes singular or ill-conditioned. In our tests we set the water level for the coefficient to $10^{-6}$ and the maximum may locally reach a value of 1 - 2, which implies a maximum-to-minimum ratio of $\approx 10^6$ for the coefficient. Despite operating at this degree of stiffness, we are able to extract accurate solutions.

{We consider flow through a medium that lies above the percolation threshold. This implies that at least one permeable channel, albeit tortuous, forms an unbroken path from one end of the domain to the other.
The percolation channel is built using a sequence of connected, permeable, smaller pieces which are parametrised according to their length, length-to-width ratio and magnitudes of the $xx, xy,$ and $yy$ component permeabilities. The dimensions of these pieces are randomly chosen, within allowed limits, ensuring that the length of any individual piece be at most a tenth of the total horizontal side of the domain but also sufficiently large in order that it be accurately resolved by the numerical scheme.} The fourth-order accurate scheme we employ here loses fidelity at a third of the Nyquist sampling rate \citep[the von Neumann analysis suggests as much; also see, e.g.][]{lele92},  implying that we must resolve each feature using at least 6 grid points. This sets a minimum size for features, below which solutions are untrustworthy. 

\begin{figure}
 \includegraphics[width=\linewidth]{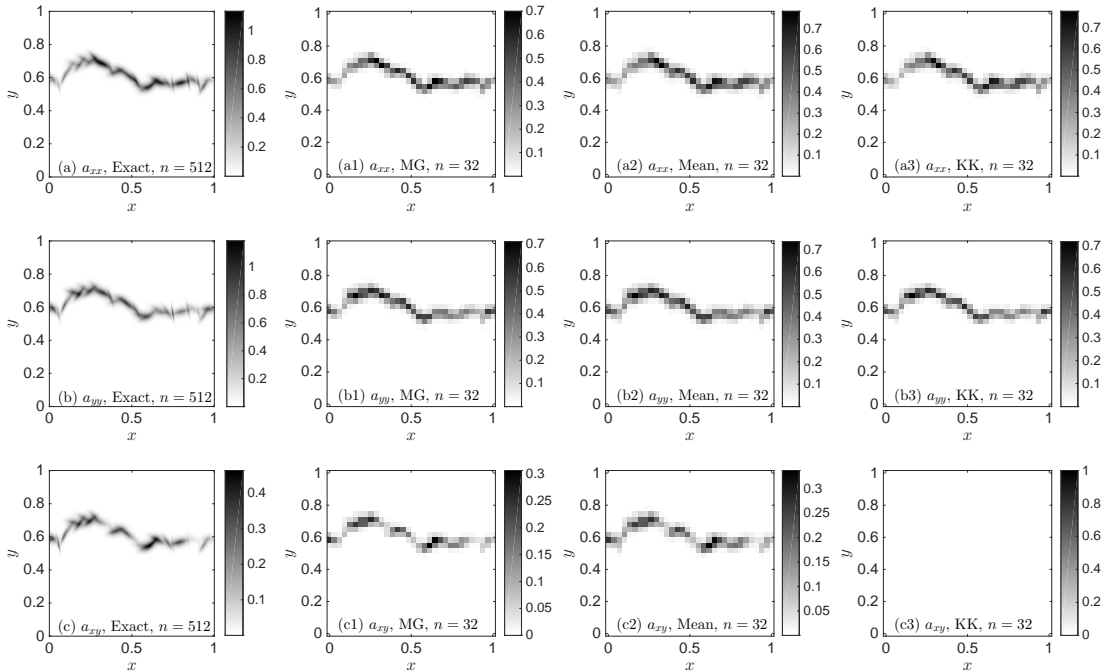}
 \caption{A single tortuous channel is seen to cross from one horizontal end of the computational box to the other of a model coefficient that has crossed the percolation threshold. The left-most column (panels a, b, c) shows the ``real" coefficient (resolution of $512\times512$), and the three successive columns show coefficients decimated using MG (a1, b1, c1), Mean (a2, b2, c2) and KK (a3, b3, c3) methods to a size of $32\times32$. Note that the colour bars of the different panels are all different and that, with successive upscaling, the magnitudes of the coefficients continue to decrease. The background permeability for the $a_{xx}$ and $a_{yy}$ coefficients is $10^{-6}$ whereas the background for the $a_{xy}$ coefficient is zero. The flow associated with this model is shown in Figure~\ref{solution_example}.\label{model.egs}}
 \end{figure}

\begin{figure}
 \includegraphics[width=0.8\linewidth]{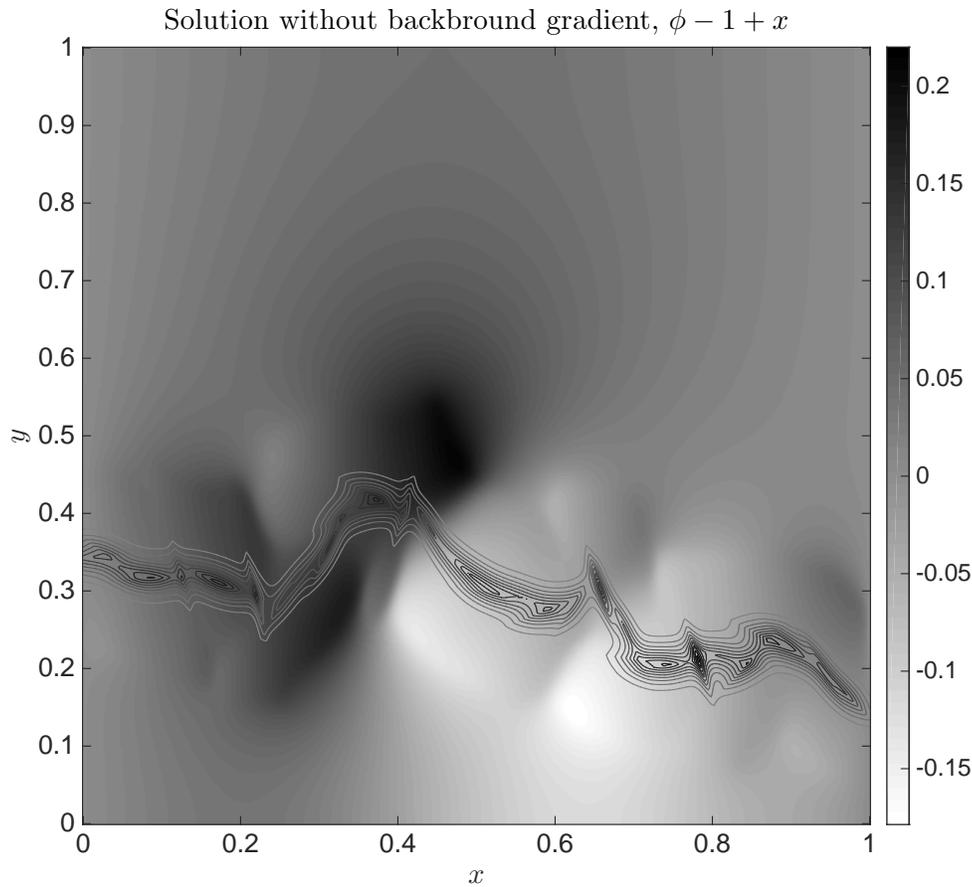}
 \caption{The flow pattern associated with applying a pressure gradient in the $x$ direction, i.e. $\phi(x=0, y) = 1$, $\phi(x=1, y) = 0$, and impermeability conditions on the upper and lower boundaries, $\partial\phi/\partial y |_ {(x,y=0)} = 0 = \partial\phi/\partial y |_ {(x,y=1)}$. The over-plotted contours show the ``exact'' coefficient $a_{xx}$ from the top-left panel of Figure~\ref{model.egs} to assist in identifying flow features. To obtain $\phi$, we solve equation~(\ref{full.eq}) using PETSc \citep{petsc-main}. Nominally, for constant permeability, we have $\phi =1-x$. In order to highlight the spatial complexity of the solution, we display $\phi - 1 + x$. { The choices that determine the size and orientation of the percolation channel (randomly assigned within allowed ranges) have an impact the total flow rate. To avoid being overly sensitive to the associated finite-size effects, we perform a statistical survey over a broad range in sizes and orientations.} \label{solution_example}}
 \end{figure}

{ The inclination of each piece with respect to the axes and coefficient magnitude is also randomly assigned, within a fixed range - inclinations must be chosen so that the percolation channel makes effective progress across the domain. The finite dimensions of these pieces can directly impact the overall flow rate.  Therefore, to mitigate the sensitivity to any individual choices in these permeability models, we generate 500 models of permeability coefficients, each with at least one percolating channel, which is in turn built using randomly chosen channel paths, sizes and permeability magnitudes.} We show an examples of one such coefficient model above the percolation threshold in Figure~\ref{model.egs} and the corresponding solution in Figure~\ref{solution_example}.

\subsection{Upscaling methods and cost}\label{upscale.method}
A number of different strategies may be employed in the process of upscaling \citep[e.g. see][and references therein]{karim}. { Coarsening a highly resolved coefficient first requires a target resolution, e.g. we might wish to upscale a $512\times512$ model to $32\times32$. Secondly, this coarsening may be accomplished either all at once or through a sequence of steps. For instance, we might decimate an $N \times N$-sized model of a coefficient by operating on blocks of size $N_{\rm block}\times N_{\rm block}$, reducing each lock to an effective coefficient. This would spatially coarsen the model of the coefficient to size $N/N_{\rm block}\times N/N_{\rm block}$ (assuming that $N_{\rm block}$ is a divisor of $N$). The process is repeated $K$ times on successively spatially reduced coefficients, where $K =\log_{N_{\rm block}}(N/{N_{\rm target}})$, eventually arriving at the target grid size $N_{\rm target}\times N_{\rm target}$ (note that the goal is $N_{\rm target} \ll N$). In 2-D, one upscaling sweep with MG requires inverting $(N / N_{\rm block})^2$ matrices of size $N^2_{\rm block} \times N^2_{\rm block}$. The cost of each inversion scales as $\left(N^2_{\rm block}\right)^3$ and the first sweep (i.e. going from $N\times N$ to $N/N_{\rm block}\times N/N_{\rm block}$), which requires $(N / N_{\rm block})^2$ inversions, has an associated computational cost that scales as $(N / N_{\rm block})^2 \times N^6_{\rm block} = N^2 N^4_{\rm block}$. If the total number of required sweeps $K$ is larger than 1, then the computational cost associated with the subsequent sweep (going from $N/N_{\rm block}\times N/N_{\rm block}$ to $N/N^2_{\rm block}\times N/N^2_{\rm block}$) is $N^2\,N^2_{\rm block}$ and so on.   
Assuming that the block size $N_{\rm block}$ is the same at each sweep (note that this is not necessary), the total cost is $N^2 N^4_{\rm block}\sum_{i=0}^{K-1} N^{-2k}_{\rm block}$. 

Consider the case of $N=512$, $N_{\rm target}  = 32$ and the simplest possible block size of $N_{\rm block} = 2$.  Assuming that inverting an $N^2_{\rm block} \times N^2_{\rm block} = 4\times4$ matrix can be accomplished in $10\mu$s, sequentially upscaling a $512\times512$ model to $32\times32$ requires 4 sweeps (from 512 to 256, 256 to 128, 128 to 64 and finally 64 to 32), and typically takes $0.5$s or so in \texttt{MATLAB}. To illustrate the escalating cost with block size, consider $16\times16$ blocks, i.e.  $N^2_{\rm block} = 256$. One sweep over the $512\times512$ grid with this block size will reduce it to $32\times32$ and the computational time, using the scaling described above, would take some $512$s. Therefore, although the number of sweeps required decreases with increasing block size, the cost of the matrix inversion dramatically increases, resulting in potentially long compute times. 

Our empirical tests (fortuitously) revealed that in the limit of small block size, i.e. $N_{\rm block} \ll N$ (we tried $N_{\rm block} = 2, 4, 8, 16$ for $N = 512$), the eventual accuracy in retrieving the flow rate does not change very much with block size. We may thus employ the computationally cheap $N_{\rm block} = 2$ strategy. 

It is important to keep in mind that computing solutions to Darcy's equation comes down to solving an implicit boundary-value problem, which translates to inverting a matrix of size $N^2\times N^2$ for a grid size $N\times N$. The computational cost of a direct matrix inversion scales as $O((N^2)^3) = O(N^6)$, although indirect solves using sophisticated sparse-matrix techniques such as multigrid can significantly mitigate the cost. For instance, solutions on a $512\times512$ grid take about 60 seconds, whereas the solution on a $32\times32$ grid takes just a few miliseconds. While $60$ seconds can still be considered as reasonably fast, the time required to solve for the flow through a permeability with sharp contrasts in 3-D on a $1024^3$ is a factor of about $32^9 \approx 10^{13}$ more expensive than a solution on a $32^3$ grid. Running MG to upscale the coefficient is therefore a trivial cost to pay to obtain accurate solutions, especially when using the lowest block size. We study three upscaling techniques here, } 
\begin{enumerate}
\item {MG: codes, written in {\texttt{MATLAB}}, take as input $N\times N$-sized (grid) tensor permeabilities and decimate them down to desired $N_{\rm target} \times N_{\rm target}$-sized coefficients. The MG analysis in Section~\ref{RG.theory} relies on the Fourier basis; however, one may equally derive the technique using the Haar wavelet basis. The computational cost scales as $N^2 N^4_{\rm block}\sum_{i=0}^{K-1} N^{-2k}_{\rm block}$, where $K = \log_{N_{\rm block}} N$.
\item Straight mean: we replace blocks of size $N_{\rm block}\times N_{\rm block}$ by arithmetic averages for each of the component coefficients $xx, xy,$ and $yy$. Thus, one level of decimation causes the size of the coefficient to drop from $N\times N$ to $N/N_{\rm block} \times N/N_{\rm block}$. Depending on the target grid size, several sweeps may need to be performed. It may be shown that the total cost scales as $N^2 N_{\rm block}^2 \sum_{k=0}^{K-1} N^{-2k}_{\rm block}$, which evidently requires fewer operations, notwithstanding a proportionality constant, than MG. For sufficiently small block size $N_{\rm block}$, however, there is little practical difference in cost between the methods.
\item KK method: is an algorithm designed to replace 2$\times$2-sized blocks of $xx$ and $yy$ components by combinations of their arithmetic and geometric means (see \citep{karim} and appendix~\ref{RG.k2}). { By construction, KK can only handle a block size of $N_{\rm block}=2$, i.e. KK is designed specifically to decimate $2\times2$-sized blocks $(a_{xx}, a_{yy})$, at adjacent $x$- and $y$-grid points, to corresponding effective values. The algorithm does not account for $xy$ coefficients, i.e. $a_{xy} = 0$ pre- and post upscaling.} For sake of completeness and comparison, we state the analytical expressions for KK and corresponding $k=1$ decimation using MG in appendix~\ref{RG.k2}. The computational cost is similar to that of the straight mean.}
\end{enumerate}

\section{Statistical study and Discussion}
It is to be expected that the less the coefficient is upscaled, the greater the accuracy of the corresponding solution. Owing to computational restrictions, the maximum-allowed size for the coarsened coefficient is set at $32\times32$ (user defined), and we therefore operate at this limit. We decimate the coefficient down to this spatial size, thereby defining the extent of coarsening. 

The size of the block, i.e. the $N_{\rm block}\times N_{\rm block}$ tile that is decimated to one tensor (following the notation in Section~\ref{upscale.method}), is also an important choice. Based on a number of tests (not discussed here) for the flow rate, we find that the error is mostly insensitive to the choice of block size, at least when it is much smaller than the entire grid dimension of the coefficient, i.e. $N_{\rm block} \ll N$. In fact, preliminary calculations reveal that $N_{\rm block}=2$ is the best choice although more thorough testing may be required to support this. Setting $N_{\rm block}=2$, we will need to decimate the $512\times512$-sized model four times to arrive at the eventual grid size of $32\times32$.

We must choose the basis on which to project the coefficients; based on  preliminary tests, we proceed with the Fourier basis, finding it to be more accurate (in terms of flow-rate prediction) than the Haar. Indeed, more careful testing is required to justify this choice.

{We perform simulations for 500 models of permeabilities, each possessing a randomly generated tortuous channel from the left boundary to the right, ensuring that we have a percolation pathway running through the medium.} For each realization, we compute the ``exact" solution for the original $512\times512$ model, decimate the model according to MG, KK and the straight mean prescriptions to various levels, i.e. $256\times256$, $128\times128$ and so on. Because of the sparsity of permeable channels, decimating below a target grid size $N_{\rm target}\times N_{\rm target} = 32\times32$ causes the information about the channel to be lost. { To avoid confusion, we recall here that the block size controles the number of steps in the upscaling process, whereas the target size denotes the extent to which upscaling is performed. Too much upscaling (i.e. when the target grid size is too small) causes the permeability coefficient to fall to the background value ($10^{-6}$ in this case), and the medium is deemed impermeable. As for the block size, small $N_{\rm block}$ implies low computational cost and vice versa.}

 \begin{figure}
 \includegraphics[width=\linewidth]{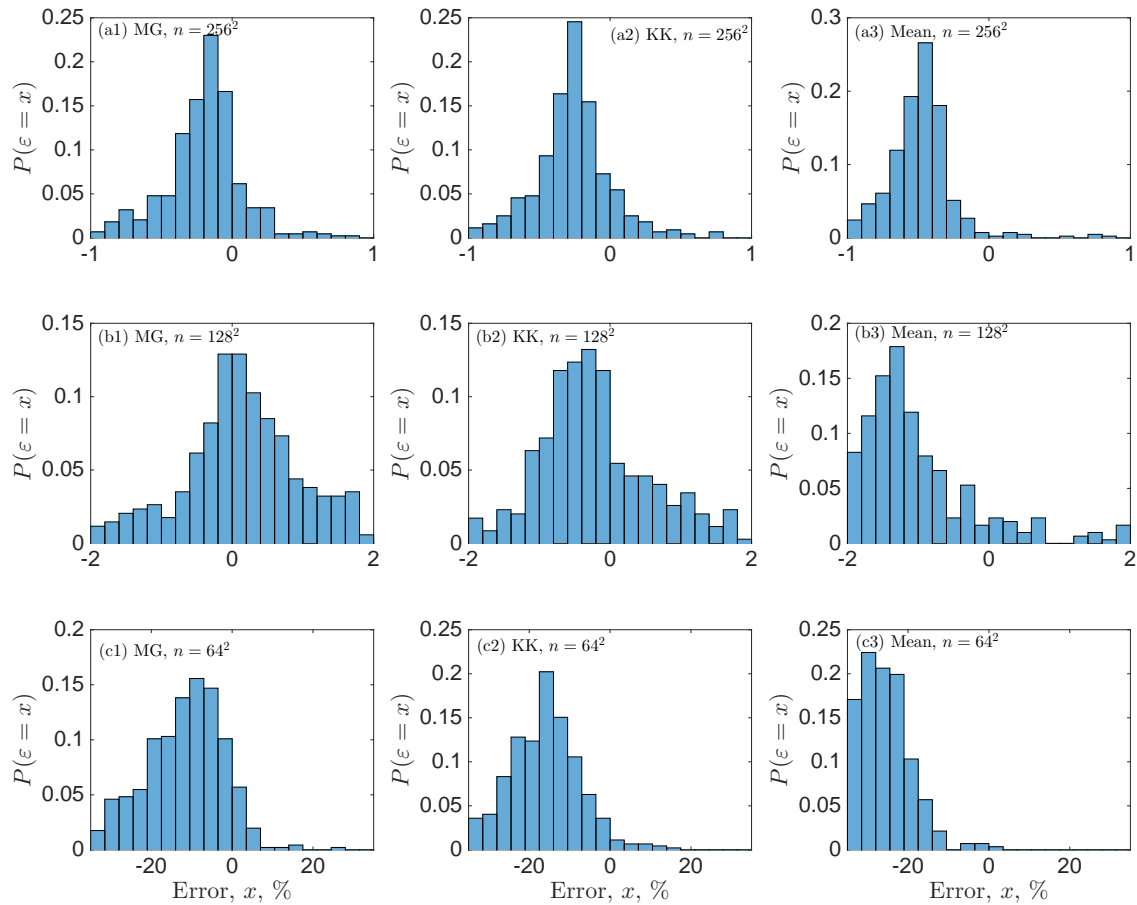}
 \caption{The statistical distribution of numerical errors associated with upscaling 500 models of permeability using three techniques, MG (panels a1, b1, c1), KK (panels a2, b2, c2) or straightforward averaging (panels a3, b3, c3). 
 In these tests, we explicitly set $a_{xy} = 0$ at the highest resolution ($512\times512$). The $x$ axis is error in predicting the correct flow rate (Eq.~[\ref{error.def}]) and the $y$ axis is the probability of incurring that error. The original coefficient is sampled at $512\times512$ resolution and subsequently decimated to resolutions of $256\times256$ (row a), $128\times128$ (row b) and finally $32\times32$ (row c). With greater rates of decimation, MG introduces a non-trivial $a_{xy}$ coefficient, whereas KK and straight mean do not address this aspect. The models lie above the percolation threshold, with a single channel linking the two horizontal boundaries. The channel may assume a tortuous path. It is seen that although all the methods show a systematic theoretical bias towards overpredicting the flowrate, MG is the best performing method\label{flow_errors_zero}}
 \end{figure}

 \begin{figure}
 \includegraphics[width=\linewidth]{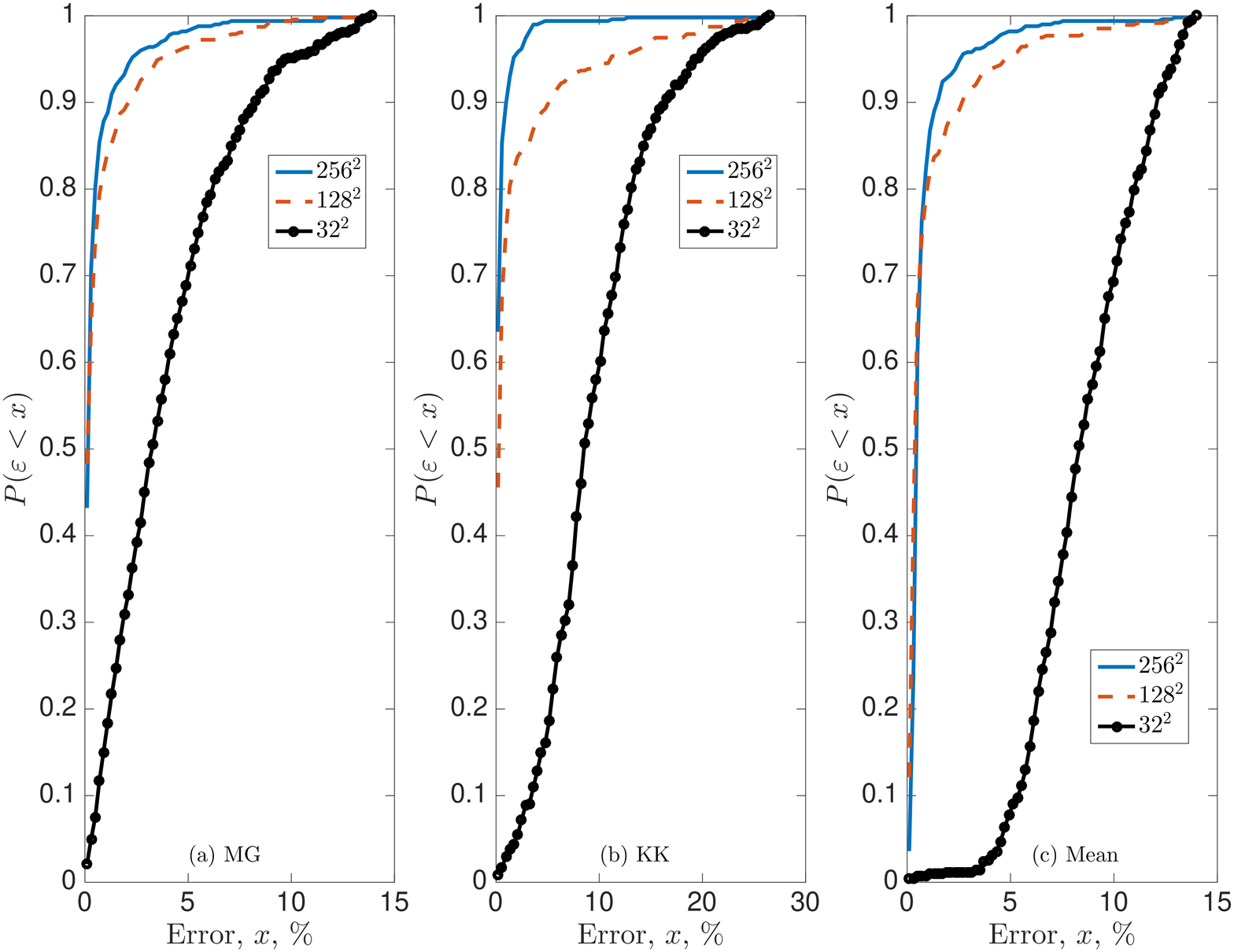}
 \caption{Probability of incurring an error less than $x$ in calculating the flow rate (Eq.~[\ref{error.def}]) when decimating from $512\times512$ to various resolutions for the three methods being considered here, MG (panel a), KK (panel b) and straight mean (panel c). This is concluded from a statistical survey of 500 randomly generated models that are above the percolation threshold. The off-diagonal coefficient $a_{xy}=0$ for the original model at the highest resolution of $512\times512$. The upscaled coefficients obtained using MG have finite $a_{xy}$ and are therefore of higher accuracy than KK and the straight mean, which do not take these into account.\label{methods-cumulative-zero}}
 \end{figure}

 \begin{figure}
 \includegraphics[width=\linewidth]{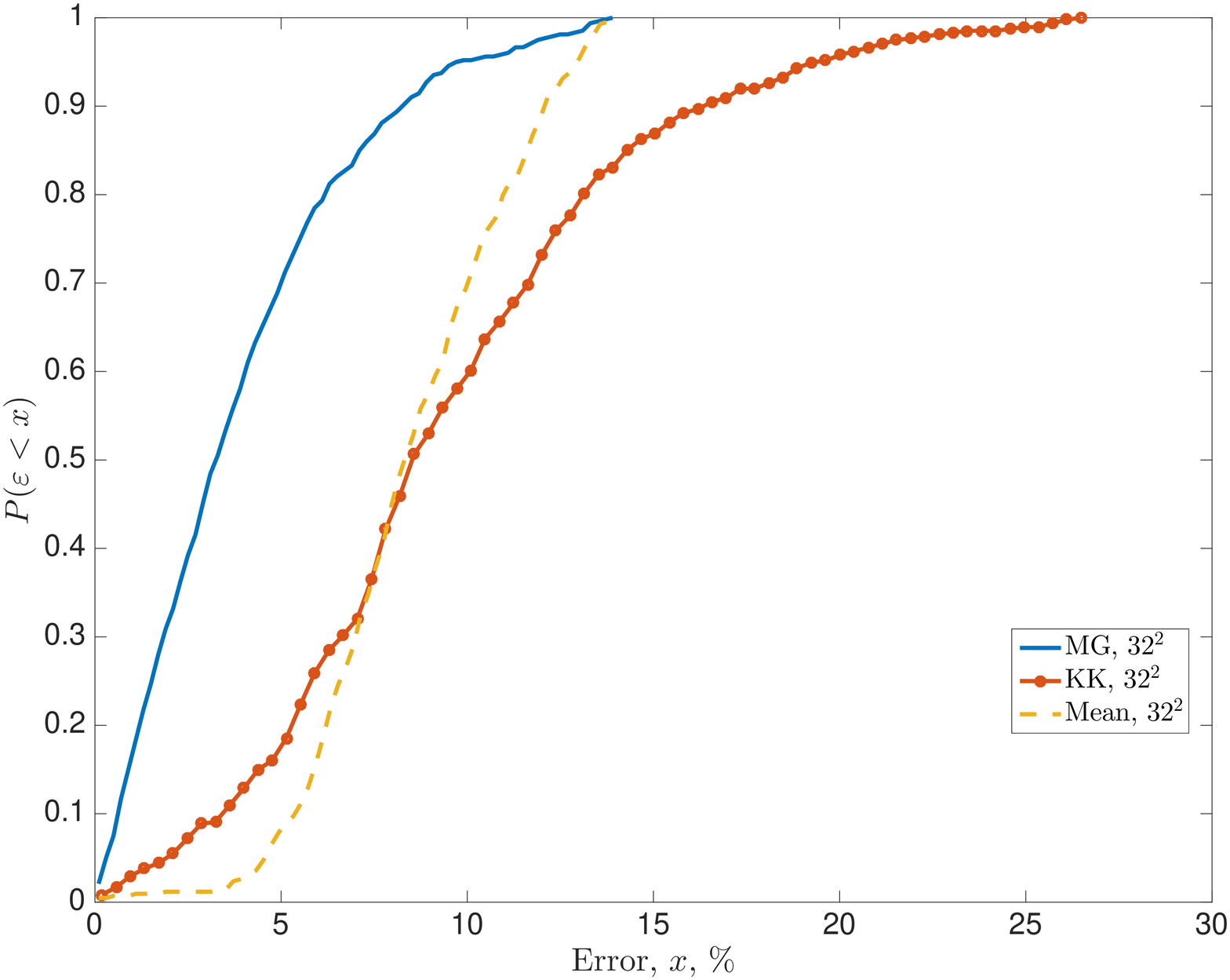}
 \caption{Probability of incurring an error less than $x$ in calculating the flow rate (Eq.~[\ref{error.def}]) when decimating from $512\times512$ to $32\times32$ for the three methods being considered here, MG, KK and straight mean. This is concluded from a statistical survey of 500 randomly generated models that are above the percolation threshold. The off-diagonal coefficient $a_{xy} = 0$ for the original model at the highest resolution of $512\times512$. When decimating down to lower resolutions, only MG is able to address the tensorial nature of the coefficient and thereby introduces finite $a_{xy}$ at $32\times32$, resulting in greater accuracy. \label{cumulative-zero}}
 \end{figure}

 \begin{figure}
 \includegraphics[width=\linewidth]{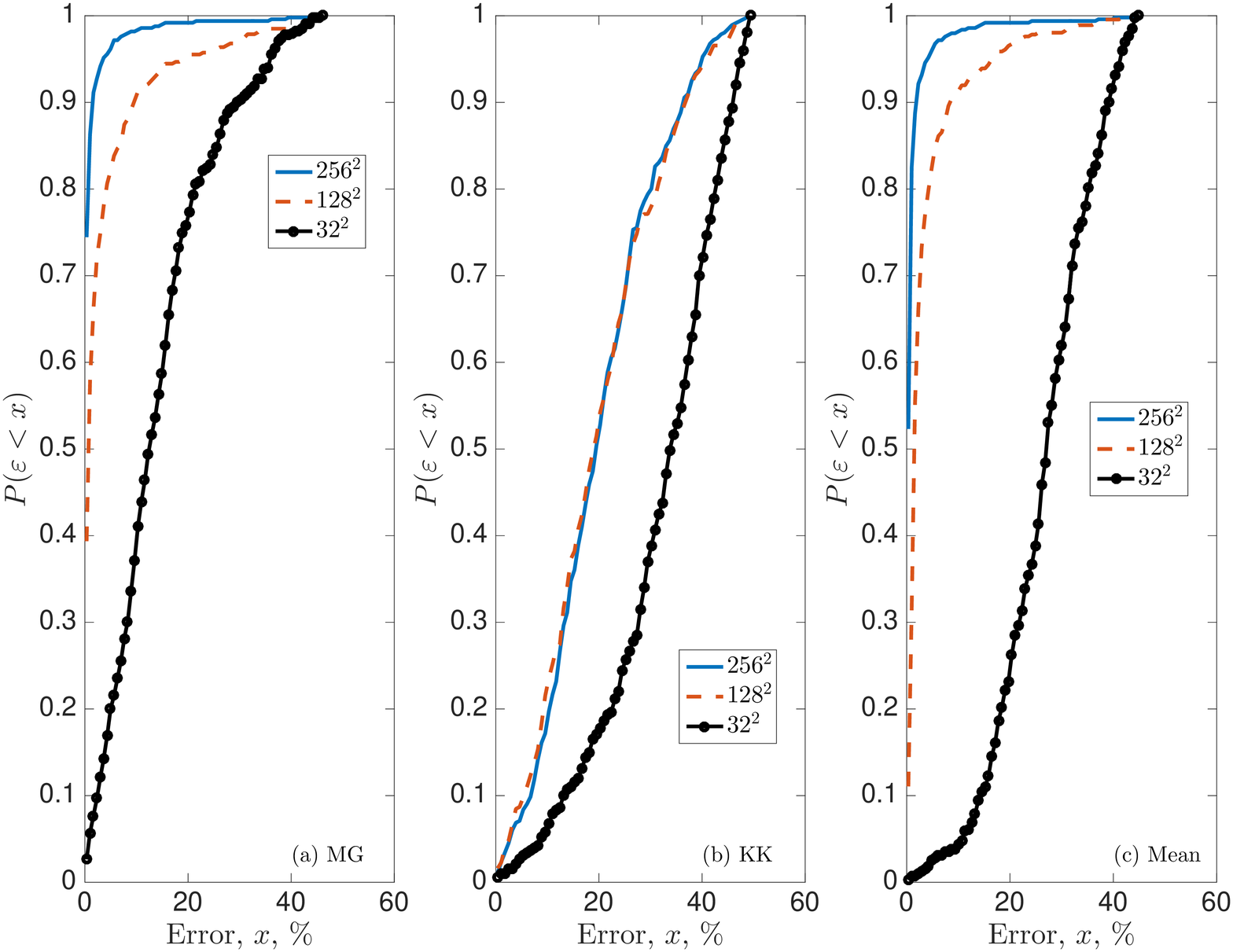}
 \caption{Probability of incurring an error less than $x$ in calculating the flow rate (Eq.~[\ref{error.def}]) when decimating from $512\times512$ to $32\times32$ for the three methods being considered here, MG (panel a), KK (panel b) and straight mean (panel c).  { The difference between this and Figure~\ref{methods-cumulative-zero} is that the off-diagonal coefficient $a_{xy}$ is now non-zero for the original model at the highest resolution of $512\times512$.}  These results are concluded from a statistical survey of 500 randomly generated models that are above the percolation threshold. The upscaled coefficients obtained using MG and straight mean techniques have finite $a_{xy}$, performing better than KK, which does not take these into account and therefore prone to significant error, even at $256\times256$.  { In general, it is seen that the errors are systematically larger for all methods when finite off-diagonal permeabilities are included (compare the $x$ axes between this and Figure~\ref{methods-cumulative-zero}).}\label{methods-cumulative-finite}}
 \end{figure}

 \begin{figure}
 \includegraphics[width=\linewidth]{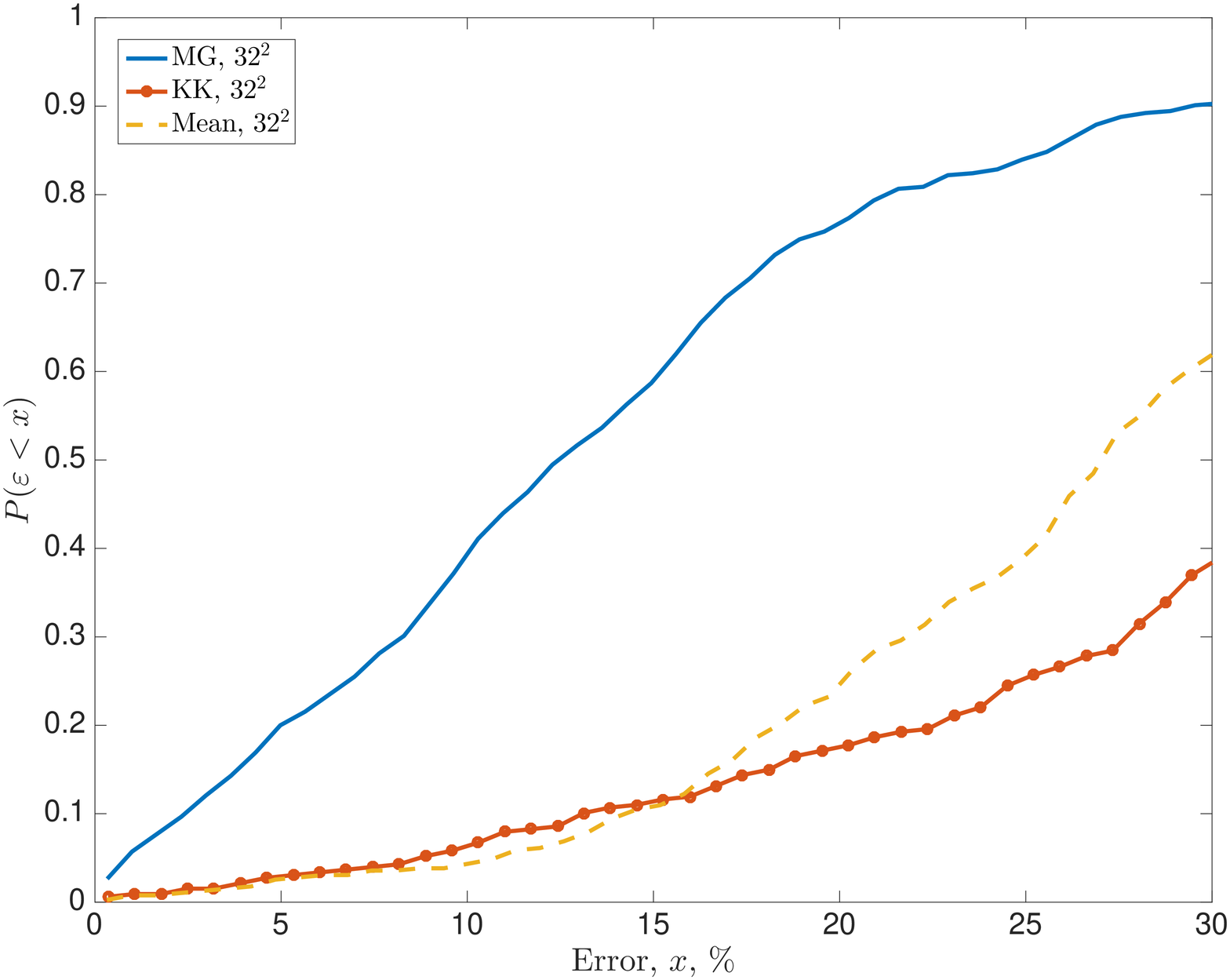}
 \caption{Probability of incurring an error less than $x$ in calculating the flow rate (Eq.~[\ref{error.def}]) when decimating from $512\times512$ to $32\times32$ for the three methods being considered here, MG, KK and straight mean. { The difference between this and Figure~\ref{cumulative-zero} is that the off-diagonal coefficient $a_{xy}$ is non-zero for the original model at the highest resolution of $512\times512$.}  The results are concluded from a statistical survey of 500 randomly generated models that are above the percolation threshold. The upscaled coefficients obtained using MG and straight mean techniques have finite $a_{xy}$ and therefore of higher accuracy than KK, which does not take these into account. { In general, it is seen that the errors are systematically larger for all upscaling methods when finite off-diagonal permeabilities are included (compare with Figure~\ref{cumulative-zero}).} \label{cumulative-finite}} 
 \end{figure}

We perform two statistical surveys, one where the the highest-resolution starting model ($512\times512$) has no off-diagonal permeability coefficient, i.e. $a_{xy}=0$ and another where it is finite. The former reflects the contemporary practice of ignoring off-diagonal coefficients in flow through porous media whereas the latter is an investigation of how such terms might influence the eventual flow estimate. { In line with expectation, we find that the accuracy in estimating the flow rates diminishes with the extent of upscaling. Note that `mild' upscaling refers to moderate $N_{\rm target} / N$ ratios, whereas `aggressive' upscaling involves pushing $N_{\rm target}$ towards smaller and smaller values. The ideal upscaling method would be very aggressive and return one effective coefficient for the entire medium, which would then describe the flow rate (see appendix~\ref{constcoeff}), thereby requiring no further computation.    
However, we find that, irrespective of the technique, upscaling to resolutions $N_{\rm target}\times N_{\rm target} = 16\times16$ and smaller results in total loss of information about the percolation channel. The effective flow rate for such highly coarsened models is comparable to the background rate of $10^{-6}$ and therefore deemed to be impermeable. We therefore halt upscaling at $32\times32$. 

Figure~\ref{flow_errors_zero} shows a histogram of probabilities of incurring a given error in the case where $a_{xy} = 0$ in the highest resolution case, as a function of the degree of upscaling. MG is seen to outperform KK slightly and the straight-mean technique significantly. This is because MG takes into account the interplay between different tensor components of the permeability. The influence between $xx, xy$ and $yy$ coefficients becomes more significant with the degree of upscaling and MG therefore is more accurate than the other two techniques. A curious feature of the coarsened solution is that it systematically overestimates the flow rate when compared to the true solution, as seen in Figure~\ref{flow_errors_zero}.

Figure~\ref{methods-cumulative-zero} plots the probability of incurring an error $\varepsilon < x$ as a function of $x$, for different methods and various extents of upscaling with $a_{xy} = 0$ at the highest resolution. Because of this, KK continues to perform reasonably well at successively lower resolutions.  Figure~\ref{cumulative-zero} tracks the probability of incurring an error $\varepsilon < x$ as a function of $x$, where the error is defined in equation~(\ref{error.def}) for upscaling to a size of $32\times32$ and with $a_{xy}=0$ at the highest resolution. MG outperforms KK and straight mean. 

Figure~\ref{methods-cumulative-finite} plots the probability of incurring an error $\varepsilon < x$ as a function of $x$, for different methods and various extents of upscaling with finite $a_{xy}$ at the highest resolution. Because KK does not take into account off-diagonal permeability terms, the degradation of its performance is immediate, evidenced by attendant errors at $256\times256$. The straight-mean method does reasonably well because successively coarsened models have finite $a_{xy}$ terms. MG continues to maintain good performance. Figure~\ref{cumulative-finite} tracks the probability of incurring an error $\varepsilon < x$ as a function of $x$, where the error is defined in equation~(\ref{error.def}), for upscaling to a size of $32\times32$ with finite $a_{xy}$ at the highest resolution. MG outperforms straightforward averaging by a wide margin (KK does poorly because it does not account for off-diagonal permeabilities). }

\section{Summary and Conclusions}
We have investigated the ability of upscaling approaches in reproducing the overall flow-rate in realistic porous media. We assembled several important elements to enable this study, such as the development of an accurate numerical solver for Darcy's equation, parametrising the permeability and building realistic models of coefficients, implementation of upscaling schemes and a statistical survey over an ensemble of coefficient models.

The analysis shows that Model-Elimination Renormalization-Group theory (MG) outperforms other upscaling techniques and the impact is higher for permeability anisotropies that are not aligned with the computational grid. One contributing reason is that MG accounts for non-diagonal terms in the permeability tensor, i.e. $a_{xy}$, an issue that gains increasing significance with greater extents of upscaling. The KK technique has no means of addressing these coefficients. Taking the straightforward mean is the coarsest method of all since no ideas pertaining to the physics or mathematics of flow through porous media are incorporated. We also find that all the techniques tend to be systematically biased towards overpredicting the flow rate. 

{ The computational time for reducing a $512\times512$ model to a $32 \times 32$ upscaled equivalent by repeatedly decimating $N_{\rm block}\times N_{\rm block} = 2\times2$ blocks using MG takes less than a second in {\texttt{MATLAB}} on a conventional desktop. Solving Darcy's equation is the most expensive step in this process, taking about 60 seconds for a $512\times512$ grid. 

The greatest gains are anticipated in 3-D, where the computation to retrieve the solution on a $1024^3$ grid is about $10^{13}$ times more expensive than on an upscaled $32^3$ model. It is therefore critical to develop accurate upscaling techniques, since the solution of Darcy's equation for complex permeabilities on large grids is computationally infeasible. Our future work entails looking at this problem in 3-D and applying upscaling methodologies to a broader set of realistic permeability models.}

\appendix
\section{Explicit coefficients}\label{RG.k2}
\subsection{Mode-Elimination Renormalization Group (MG)}
For a $2\times2$ block with $xx$ coefficients denoted by $a_{11}, a_{12}, a_{21}, a_{22}$, the $yy$ coefficients by $b_{11}, b_{12}, b_{21}, b_{22}$ and assuming the $xy, yx$ coefficients to be zero, the MG permeabilities ${\bar a}_{xx}, {\bar a}_{xy}, {\bar a}_{yy}$ are given by
\begin{eqnarray}
D &=& (a_{11} + a_{12}) (a_{21} + a_{22})(b_{11} + b_{12} + b_{21} + b_{22}) + (a_{11} + a_{12} + a_{21} + a_{22})(b_{11} + b_{21}) (b_{12} + b_{22}),\nonumber\\
N_1 &=& (a_{11}a_{12}a_{21} + a_{11}a_{12}a_{22} + a_{11}a_{21}a_{22} + a_{12}a_{21}a_{22})(b_{11} + b_{12} + b_{21} + b_{22})\nonumber\\
&+& (a_{11} + a_{21}) (a_{12} + a_{22})(b_{11} + b_{21})(b_{12} + b_{22}), \\
N_2 &=& (a_{11} + a_{12}) (a_{21} + a_{22})(b_{11} + b_{12})(b_{21} + b_{22})\nonumber\\
&+& (b_{11}b_{12}b_{21} + b_{11}b_{12}b_{22} + b_{11}b_{21}b_{22} + b_{12}b_{21}b_{22})(a_{11} + a_{12} + a_{21} + a_{22}), \\
{\bar a}_{xx} &=& \frac{N_1}{D}, \,\,\,\,\,\,\,\,{\bar a}_{yy} = \frac{N_2}{D},\\
{\bar a}_{xy} &=& - \frac{(a_{11} a_{22} - a_{12} a_{21})(b_{11} b_{22} - b_{12} b_{21})}{D} = {\bar a}_{yx}.
\end{eqnarray}

\subsection{The KK method}
The equations denoting effective conductivity are described here for the sake of completeness, taken from equations (13) and (14) of \cite{karim},
\begin{eqnarray}
{\bar a}_{xx} &=& \left[\frac{(a_{11} + a_{21})(a_{12} + a_{22})}{(a_{11} + a_{12})(a_{21} + a_{22})} \frac{a_{11} + a_{12} + a_{21} + a_{22}}{a_{11} a_{12}(a_{21} + a_{22}) + a_{21} a_{22} (a_{11} + a_{12})} \right]^{\frac{1}{2}},\\
{\bar a}_{yy} &=& \left[\frac{(b_{11} + b_{12})(b_{21} + b_{22})}{(b_{11} + b_{21})(b_{12} + b_{22})} \frac{b_{11} + b_{12} + b_{21} + b_{22}}{b_{11} b_{21}(b_{12} + b_{22}) + b_{12} b_{22} (b_{11} + a_{21})} \right]^{\frac{1}{2}},\\
{\bar a}_{xy} &=&  0 =  {\bar a}_{yx}.
\end{eqnarray}

\section{Solution for constant coefficient case}\label{constcoeff}
Given a constant permeability coefficient, the symmetry of the problem implies that there is no variation as a function $y$ (the boundary conditions on the upper and lower boundaries are both zero-Neumann and the horizontal boundaries are constant). The $y-$derivatives therefore drop out and we are left with $a_{xx} \partial^2_x \phi(x) = 0$. Given the boundary conditions of $\phi(x=0,y) = 1$ and $\phi(x=1,y) = 0$, the solution is $\phi(x) = 1 -x$. The flow rate is therefore $f = -\int_0^1 dy\, a_{xx} \partial_x\phi = a_{xx}$.

\begin{acknowledgments}

SMH acknowledges support from consultancy project PT54324 with Shell India, Ramanujan fellowship SB/S2/RJN-73/2013, and the Max-Planck Partner Group Program.
\end{acknowledgments}
\bibliography{ms.bbl}
\end{document}